\documentclass[preprint,12pt,authoryear]{elsarticle}

\usepackage{amssymb}
\usepackage{amsmath}

\usepackage{amsmath}    
\usepackage{amsfonts}   
\usepackage{amssymb}    
\usepackage{bm}         
\usepackage{algorithm}
\usepackage{algpseudocode}
\usepackage{amsmath, amssymb, graphicx}  
\usepackage{booktabs}  
\usepackage{subcaption}
\usepackage{hyperref}

\journal{}

\begin{document}

\begin{frontmatter}

\title{Adaptive Nonlinear Elimination Preconditioning for Transport in Fractured Porous Media} 

\author[1]{Omar Chaabi}
\author[1,2]{Mohammed Al Kobaisi\corref{3}}
\cortext[3]{M.S.K.AlKobaisi@tudelft.nl; mohammed.alkobaisi@ku.ac.ae}

\affiliation[1]{Khalifa University, Chemical and Petroleum Engineering, Abu Dhabi, UAE}
\affiliation[2]{Delft University of Technology, Delft Institute of Applied Mathematics, The Netherlands}

\begin{abstract}

Sequential implicit (SI) formulations are gaining increasing interest due to their ability to decouple reservoir simulation problems into distinct flow and transport subproblems, allowing for the use of specialized solvers tailored to each. This separation often improves solver efficiency and flexibility, especially in weakly coupled systems. However, for fractured reservoirs, even the decoupled subproblems may generate nonlinearly stiff systems. This is specifically evident in the transport subproblem, where fracture-induced nonlinearity imbalances often lead to poor Newton convergence, including failed iterations and frequent timestep cuts. To address this challenge, we propose and investigate an adaptive Nonlinear Elimination (NE) preconditioned exact Newton algorithm specifically tailored to transport subproblems that arise from the sequential splitting of two-phase flow in fractured porous media. The proposed method is evaluated through a series of waterflooding test cases involving discrete fracture networks. The adaptive NE-preconditioned algorithm consistently demonstrates improved convergence behavior and computational efficiency compared to standard Newton. 
\end{abstract}




\end{frontmatter}

\section{Introduction}
\label{sec1}

Reservoir simulation is vital for modern reservoir management and serves as a tool that helps engineers plan and improve subsurface energy strategies. By creating detailed digital models of underground formations, informed decisions can be made about full field development and recovery methods. However, achieving detailed accuracy often involves high-fidelity models with millions of grid cells, which can significantly increase computational demands and processing time. This high computational load makes the efficiency of simulation software for tasks such as history matching, uncertainty assessments, and forecasting essential and quite challenging. 

Of the various subsurface formations, fractured reservoirs bring additional complexities due to the extreme contrast in flow properties between fractures and matrix rock. Modeling fluid flow accurately in these systems requires fine-resolution grids to represent detailed fracture networks and their interactions with the surrounding rock. Consequently, the computational effort increases significantly, highlighting the critical need for efficient, highly performant simulation tools designed specifically for fractured reservoirs.

Traditionally, the fully implicit (FI) method has been the standard temporal discretization scheme of choice in reservoir simulation \citep{aziz1979petroleum}. The FI formulation handles all physical processes within a unified implicit framework, ensuring numerical stability albeit at the cost of solving large systems of nonlinear equations at each time step. Typically, Newton’s method is employed to handle these nonlinear systems. However, because the FI approach solves for all variables simultaneously, it inherently runs the risk of producing nonlinearly stiff systems \citep{lie2024enhancing}. That is, systems with localized nonlinearities in time and space. Such stiffness can lead to wasteful iterations, frequent timestep cuts, or the need for stringent timestep constraints to achieve convergence.

An alternative temporal discretization scheme that has been gaining more attention is the sequential implicit (SI) scheme, where the typical reservoir simulation governing equations are decoupled into flow and transport subproblems. This can be advantageous because we can assign specialized solvers targeted for each subproblem, thus adding flexibility in the choice of solution strategy and time stepping for each of the solvers \citep{kozlova2016real}.  

The early success of the multiscale finite volume method, first introduced by \citet{jenny2003multi}, played a pivotal role in motivating the development of the first sequential fully implicit (SFI) method proposed by \citet{jenny2006adaptive}. In turn, the emergence and advancements of SFI formulations spurred further research into multiscale finite volume methods \citep{hajibeygi2008iterative, zhou2011two, wang2014algebraic, moyner2016multiscale, lie2017feature, klemetsdal2020accelerating, chaabi2024algorithmic}. Over time, the interplay between these two approaches—each addressing key challenges in reservoir simulation—led to enabling their integration into commercial simulators \citep{kozlova2016real, lie2016successful}. Given the inherent multiscale nature of fractured reservoirs, extending existing multiscale methods to such systems has been the focus of several notable efforts \citep{sandve2014physics, shah2016multiscale, bosma2017multiscale, ctene2016algebraic}, and it remains an active area of research \citep{watanabe2023strongly, kumar2024algebraic, watanabe2025physics}. To the best of our knowledge, all these efforts were primarily concerned with efficient solvers for the flow subproblem.

Unfortunately, in fractured reservoirs, adopting a Sequential Implicit (SI) scheme does not entirely eliminate the risk of encountering nonlinearly stiff systems. This is mainly due to the sharp saturation changes that can develop near fractures, driven by the stark contrast in flow properties between highly permeable fractures and the surrounding rock matrix. Because of these imbalanced nonlinearities the Newton solver may experience slow convergence in the form of failed iterations and time step cuts. Nonlinear preconditioning can reduce these issues by solving nonlinear subproblems in inner iterations or by performing inner relaxation for problematic variables.

Analogously to linear preconditioning, a nonlinear preconditioner can be applied to the left or right of a nonlinear function, depending on the desired formulation and solution strategy \citep{liu2018note}. Examples of left preconditioning include the Additive and Multiplicative Schwarz Preconditioned Inexact Newton methods, known respectively as ASPIN \citep{cai2002nonlinearly} and MSPIN \citep{liu2015field}, as well as a Restricted Additive Schwarz Preconditioned Exact Newton method called RASPEN \citep{dolean2016nonlinear}. Left preconditioning transforms the original stiff nonlinear system into a better conditioned system with more balanced nonlinearities and subsequently solves it using a Newton-like approach. Such local-global nonlinear domain-decomposition methods have been applied in studies of two-phase flow in porous media \citep{skogestad2013domain,liu2013fully,skogestad2016two}. Other notable efforts related to porous media flow include extending ASPEN to three-phase compositional flow \citep{klemetsdal2022numerical}, developing an adaptive ASPEN-based hybrid implicit solution formulation \citep{klemetsdal2022adaptive}, and proposing an effective formulation that achieves computational speedup by adaptively alternating between ASPEN and Newton methods during the simulation \citep{moyner2024nonlinear}. In contrast, right nonlinear preconditioning—such as nonlinear elimination (NE) preconditioning \citep{lanzkron1996analysis, cai2011inexact}—preserves the original nonlinear function and instead modifies the unknowns of the nonlinear system, typically by performing subdomain solves or applying inner relaxation to problematic variables. The application of NE preconditioning can be viewed as a subspace correction step that potentially provides a better starting point for the global Newton iteration. Variants of the NE preconditioner have been applied to two-phase flow in porous media \citep{yang2016active, yang2017nonlinearly, luo2021nonlinear}, three-phase flow with capillarity \citep{yang2020nonlinearly}, steady-state incompressible flow problems \citep{luo2020multilayer}, and thermal-hydraulic-mechanical (THM) simulations of fractured reservoirs \citep{wang2015efficient}.

In this paper, we address the slow computational performance of Newton solvers encountered in transport subproblems that arise from the sequential implicit splitting of two-phase flow in fractured reservoirs. Our investigation attributes this slowdown to nonlinearity imbalances caused by abrupt flux variations across fractures. To overcome this issue, we propose an adaptive NE preconditioned exact Newton algorithm tailored for two-phase flow in fractured reservoir simulations. 

\section{Mathematical Model and Governing Equations}
\label{sec2}

In the following, we consider two-phase flow in fractured reservoirs modeled using the embedded discrete fracture method \citep{lee2001hierarchical, li2008efficient, moinfar2012development}, based on the implementation available in the open-source Matlab Reservoir Simulation Toolbox \citep{lie2021advanced}. For a d-dimensional computational domain \(\Omega\) containing \( N_f \) discrete fractures, fluid flow within the matrix domain is governed by mass conservation equations indexed by \( m \), 
\begin{equation}
\frac{\partial (\phi^m \rho_\alpha S^m_\alpha)}{\partial t} + \nabla \cdot (\rho_\alpha \mathbf{v}^m_\alpha) = \rho_\alpha \left[ q^m_{\alpha} - \sum_{i=1}^{N_f} q^{m,i}_{\alpha} \right], \quad \Omega_m \subset \mathbb{R}^d.
\end{equation}
Whereas fluid flow within the fractures is governed by a system of \( N_f \) mass conservation equations. Each fracture is identified by an index \( i \in [1, N_f] \). The governing equation for fractures take the following form:
\begin{equation}
\frac{\partial (\phi^i \rho_\alpha S^i_\alpha)}{\partial t} + \nabla \cdot (\rho_\alpha \mathbf{v}^i_\alpha) = \frac{\rho_\alpha}{a_i} \left[ q^i_{\alpha} - q^{i,m}_{\alpha} - \sum_{\substack{j=1, j \neq i}}^{N_f} q^{i,j}_{\alpha} \right], \quad \Omega_i \subset \mathbb{R}^{d-1}.
\end{equation}

Here, \( \phi^m \) and \( \phi^i \) denote the porosity of the matrix and fracture, respectively. The term \( \rho_\alpha \) represents the density of fluid phase \( \alpha \). The saturation of fluid phase \( \alpha \) in the matrix and in fracture \( i \) is denoted by \( S^m_\alpha \) and \( S^i_\alpha \), respectively. The fracture aperture is given by \( a_i \). Moreover, \( \mathbf{v}^m_\alpha \) and \( \mathbf{v}^i_\alpha \) represent the Darcy velocity of fluid phase \( \alpha \) in the matrix and fracture medium \( i \), respectively. The terms \( q^i_\alpha \) and \( q^m_\alpha \) are source or sink terms in the fractures and matrix, respectively. The transfer function \( q^{im}_\alpha \) represents flow from fracture \( i \) to the matrix, while \( q^{mi}_\alpha \) represents the reverse direction—from matrix to fracture. Additionally, \( q^{ij}_\alpha \) denotes flow from fracture \( i \) to fracture \( j \). 

All fracture equations are coupled to the matrix through the terms \( q^{im}_\alpha \) and \( q^{mi}_\alpha \), while coupling between fractures is represented by \( q^{ij}_\alpha \). The velocity \( \mathbf{v}_\alpha \) of fluid phase \( \alpha \) follows Darcy’s law:

\begin{equation}
\mathbf{v}_\alpha = - \frac{k_{r\alpha}}{\mu_\alpha} \mathbf{K} \left( \nabla p_\alpha - \rho_\alpha \mathbf{g} \right)
\end{equation}
Here, \( k_{r\alpha} \) is the relative permeability of fluid phase \( \alpha \), \( \mu_\alpha \) is the viscosity of phase \( \alpha \), \( \mathbf{K} \) is the absolute permeability tensor, \( p_\alpha \) is the pressure of phase \( \alpha \), and \( \mathbf{g} \) denotes the gravitational acceleration vector.

\begin{figure}[!hbt]
    \centering
    \includegraphics[width=1\textwidth]{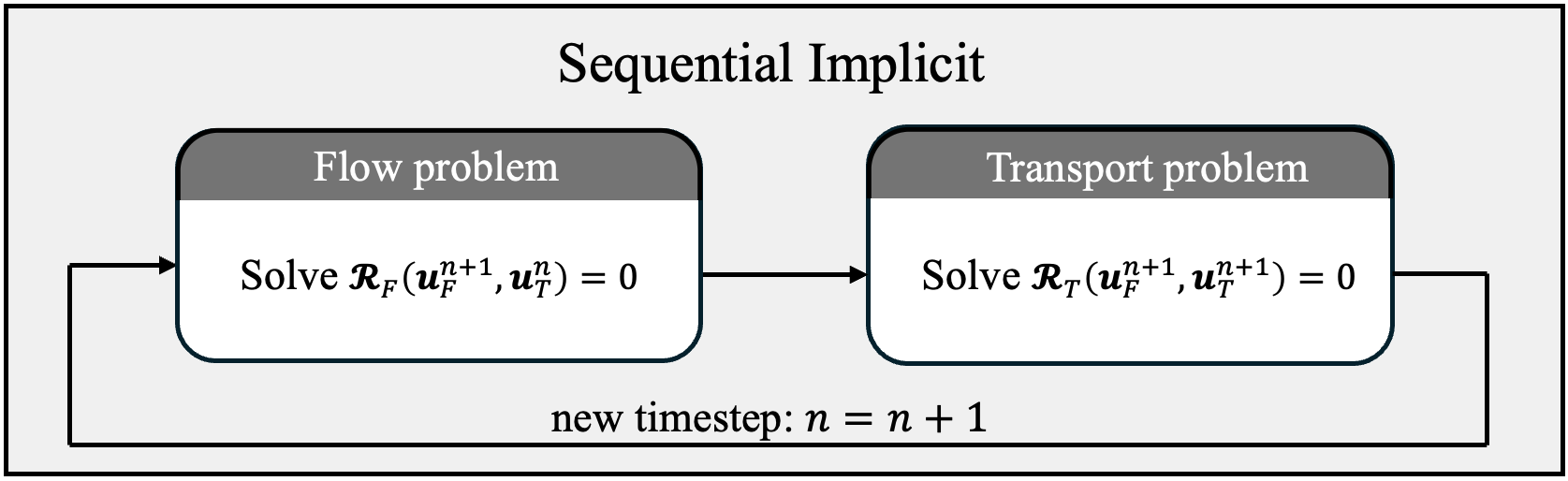}
    \caption{Schematic of the Sequential Implicit (SI) solution strategy.}
    \label{Figure 1}
\end{figure}

Having introduced the governing equations, it is necessary to define a strategy for advancing the solution in time. Various solution strategies were briefly discussed in the introduction; in this work, we focus exclusively on the sequential implicit method. For a detailed overview of alternative strategies, the reader is referred to \citet{klemetsdal2019efficient}. In SI schemes, the governing equations are decomposed into separate flow and transport subproblems. Advancing the solution in time therefore requires solving two decoupled nonlinear algebraic systems at each time step. Figure \ref{Figure 1} illustrates the sequence of steps involved in advancing the solution over a single timestep using the SI method. The figure shows our governing equations formulated in a residual-based compact vector form:
\begin{equation}
\boldsymbol{\mathcal{R}}(\mathbf{u}) = \mathbf{0}. 
\end{equation}
The vector \( \textbf{u} \) is the set of primary unknowns associated with a spatial dimension in our computational domain \( \Omega \) that depends on whether we are solving for flow \( \boldsymbol{\mathcal{R}}_F(\textbf{u}_F) = 0 \) or transport \( \boldsymbol{\mathcal{R}}_T(\textbf{u}_T) = 0 \). The standard approach is to use Newton’s method to solve each subproblem. Expanding the residual equation around an initial guess \( \bm{\textbf{u}}^0 \) using Taylor series and assuming that a corrected set of primary unknowns \( \bm{\textbf{u}}^0 + \Delta \bm{\textbf{u}} \) satisfies the residual equation, we get:

\begin{equation}
\textbf{0} = \boldsymbol{\mathcal{R}}(\bm{\textbf{u}}^0 + \Delta \bm{\textbf{u}}) = \boldsymbol{\mathcal{R}}(\bm{\textbf{u}}^0) + \textbf{J} \Delta \bm{\textbf{u}} + \mathcal{O}(\Delta \bm{\textbf{u}}^2),
\end{equation}
where \( \textbf{J} \) denotes the Jacobian matrix of the residual. Neglecting higher-order terms leads to:

\begin{equation}
\bm{\textbf{u}}^{(k+1)} = \bm{\textbf{u}}^{(k)} + \Delta \bm{\textbf{u}}, \quad \Delta \bm{\textbf{u}} = -\textbf{J}^{-1} \boldsymbol{\mathcal{R}}(\bm{\textbf{u}}^{(k)}).
\end{equation}
This equation represents the iterative scheme of a Newton-Raphson nonlinear solver with superscripts \( ^{(k+1)} \) and \( ^{(k)} \) indicating next and current iteration levels, respectively. 

Next, we present a right preconditioning technique for Newton’s method, with a specific focus on the transport subproblem \( \boldsymbol{\mathcal{R}}_T \) in fractured reservoirs. In particular, our goal is to improve the performance of the transport solver by solving a decoupled fractured transport subproblem, which serves as a preconditioner for the global Newton transport solver.

\section{Nonlinear Elimination preconditioned Exact Newton algorithm}
\label{sec3}
In this section, we present a nonlinear elimination preconditioned Exact Newton algorithm for solving transport problems in fractured reservoirs. Consider a transport nonlinear residual function \( \boldsymbol{\mathcal{R}}_T(\textbf{u}) \), where \( \textbf{u} \in \mathbb{R}^n \) is split into two subsets of variables, \( \textbf{u} = (\textbf{u}_m, \textbf{u}_f) \). An index set \( \mathcal{S} \) is likewise partitioned into two subsets: \( \mathcal{S}_m \) for the matrix and \( \mathcal{S}_f \) for the fracture, such that \( \mathcal{S} = \mathcal{S}_m \cup \mathcal{S}_f \). The nonlinear system can then be reformulated in a coupled block matrix form as follows:

\begin{equation}
\boldsymbol{\mathcal{R}}(\textbf{u}) = 
\begin{bmatrix}
\boldsymbol{\mathcal{R}}_{mm}(\textbf{u}_m, \textbf{u}_f) & \boldsymbol{\mathcal{R}}_{mf}(\textbf{u}_m, \textbf{u}_f) \\
\boldsymbol{\mathcal{R}}_{fm}(\textbf{u}_m, \textbf{u}_f) & \boldsymbol{\mathcal{R}}_{ff}(\textbf{u}_m, \textbf{u}_f)
\end{bmatrix}
\begin{bmatrix}
\textbf{u}_m \\
\textbf{u}_f
\end{bmatrix}
= \textbf{0}.
\label{Global Nonlinear system}
\end{equation}
For this partition, we define a subspace \( \mathcal{W}_f \subset \mathbb{R}^n \) corresponding to the set \( \mathcal{S}_f \) as  
\begin{equation}
\mathcal{W}_f = \left\{ \textbf{w} \mid \textbf{w} = (w_1, \ldots, w_n)^T \in \mathbb{R}^n,\; w_i = 0 \text{ if } i \notin \mathcal{S}_f \right\},
\end{equation}
and the corresponding mapping operator \( \Pi_f \), which maps a vector from \( \mathbb{R}^n \) to \( \mathcal{W}_f \). Similarly, we define a subspace  
\begin{equation}
\mathcal{W}_m = \left\{ \textbf{w} \mid \textbf{w} = (w_1, \ldots, w_n)^T \in \mathbb{R}^n,\; w_i = 0 \text{ if } i \notin \mathcal{S}_m \right\},
\end{equation}
and the corresponding operator \( \Pi_m \), which maps vectors from \( \mathbb{R}^n \) to \( \mathcal{W}_m \). Analogously to \( \Pi_f \), we define a residual restriction operator \( \Psi_f \) that acts on the full residual \(\boldsymbol{\mathcal{R}}(\textbf{u})\) and extracts  a nonlinear fracture subdomain problem. Neumann boundary conditions are imposed on \(\partial\Omega_f\) and derived from the matrix-fracture interface fluxes. By applying this restriction to our nonlinear system, eq. \ref{Global Nonlinear system}, we get:
\begin{equation}
\boldsymbol{\mathcal{R}}_{\mathcal{S}_f}(\textbf{u}) = \Psi_f\left( \boldsymbol{\mathcal{R}}(\textbf{u}) \right) = \begin{bmatrix}
\boldsymbol{\mathcal{R}}_{fm}(\textbf{u}) & \boldsymbol{\mathcal{R}}_{ff}(\textbf{u})
\end{bmatrix}.
\end{equation}
Here, we introduce a simplifying assumption, while performing local subdomain solves, by considering a completely decoupled fracture system. Under this assumption, the nonlinear residual corresponding to the fracture index set \( \mathcal{S}_f \) is approximated by isolating its local contribution. Specifically, we approximate \( \boldsymbol{\mathcal{R}}_{\mathcal{S}_f}(\textbf{u}) \) as:
\begin{equation}
    \boldsymbol{\widetilde{\mathcal{R}}}_{\mathcal{S}_f}(\textbf{u}) = \boldsymbol{\mathcal{R}}_{ff}(\textbf{u}).
\end{equation}
This approximation neglects the coupling terms between the matrix and fracture domains, and assumes that the residual in the fracture is governed primarily by fracture-fracture interactions. Such a simplification reduces computational complexity, and is generally reasonable for time-steps where the flow is predominantly taking place within the fracture network which also aligns with our aim of adaptively activating our preconditioning algorithm.

\subsection{Adaptive Preconditioning}
\label{subsection3.1}
A notable drawback of nonlinear local-global preconditioning strategies lies in the added computational overhead associated with solving local subproblems. This overhead may not be justified if local solves are applied for every global Newton iteration and the benefits of preconditioning are realized only during a limited number of iterations or timesteps. In this section, we propose a physics-based adaptive criteria with the goal of taking advantage of the benefits of nonlinear preconditioning only when needed and thus reduce the total number of subproblem solves. 

Our numerical investigations revealed that, on a per-timestep basis, there is a strong correlation between abrupt changes in water flux across fracture faces and poor convergence behavior of the nonlinear solver. To capture these dynamics, we introduce two quantities that characterize temporal variations in water fluxes across matrix-fracture interfaces.

Let \( \mathcal{F} \) denote the index set of all matrix-fracture faces. At time step \( n \), we define the difference in water flux across each face \( i \in \mathcal{F} \) as:
\[
\Delta F_i^n = F_i^n - F_i^{n-1},
\]
where \( F_i^n \) denotes the water flux across face \( i \) at time step \( n \). To quantify how rapidly the flow regime is changing, we also consider a second difference quantity:
\[
\Delta^2 F_i^n = \Delta F_i^n - \Delta F_i^{n-1},
\]
which measures the rate of change of the first difference. This quantity serves as an indicator of sudden flow transitions, such as fluids entering or exiting a fracture.

Having introduced these quantities, we define an adaptive NE preconditioning criterion based on the ratio of second- to first-order differences, averaged across all matrix-fracture faces:
\[
\displaystyle\sum_{i \in \mathcal{F}} \frac{  \left| \Delta^2 F_i^n \right| }{ \left| \Delta F_i^n \right| } > \gamma,
\]
where \( \gamma \) is a threshold that determines when the NE preconditioner should be activated. A small ratio indicates relatively stable flux behavior, whereas a high ratio suggests sharp transitions in the flow regime that may hinder convergence and thus justify preconditioning.
\subsection{Remarks on Algorithmic Choices}
\label{subsection3.2}
The algorithmic decisions underlying our exact Newton nonlinear elimination (EN-NE) preconditioning strategy, outlined in the accompanying pseudo-code (Algorithm \ref{alg:adaptive-preconditioner}), were guided by a balance between effectiveness and computational simplicity. Our goal was to demonstrate that even straightforward heuristics, guided by physics, can yield meaningful improvements when integrated into transport solvers for fractured porous media. We remark here that many other strategies can be devised.

A key choice was to designate the fracture subdomain as the ``problematic region,'' where nonlinear difficulties tend to arise due to its dominant flow behavior and sharp contrast with the surrounding rock matrix. Based on this, we considered an approximate fracture residual system \( \boldsymbol{\widetilde{\mathcal{R}}}_{\mathcal{S}_f} \) and opted to perform only a limited number of nonlinear iterations on this subproblem, rather than solving it to full convergence. This results in a reduced computational cost for local subdomain solves while still providing an improved initial guess for the global Newton iteration.

In addition, we introduced a mechanism for adaptive activation of the NE preconditioner, based on a simple scalar criterion derived from temporal changes in matrix-fracture fluxes. When this criterion is satisfied for a given timestep, the preconditioner is applied only during the first nonlinear iteration, thereby avoiding the overhead of repeated application throughout the entire solve. Based on extensive experimentation, a suitable range for the threshold parameter \(\gamma\) was found to lie between 0.1 and 0.3. Lower values (e.g., \(\gamma = 0.1\)) result in more frequent activation of the preconditioner, while higher values (e.g., \(\gamma = 0.3\)) restrict activation to instances of pronounced flux variations across fracture interfaces.

\begin{algorithm}[H]
\caption{Adaptive EN-NE algorithm for fractured transport problems}
\label{alg:adaptive-preconditioner}
\begin{algorithmic}[1]
\While{$n < N$}
    \State Solve $\boldsymbol{\mathcal{R}}_F(\mathbf{u}) = 0$ using standard Newton
    \State Construct $\boldsymbol{\mathcal{R}}_T(\mathbf{u})$, based on updated unknowns
    \While{$k = 1, 2, 3, \dots$ until convergence}
        \If{$k = 1$ and $\displaystyle \sum_{i \in \mathcal{F}} \frac{|\Delta^2 F_i^k|}{|\Delta F_i^k|} > \gamma$}
            \State Extract subproblem: $\boldsymbol{\widetilde{\mathcal{R}}}_f(\mathbf{u}_f) = \Psi_f(\boldsymbol{\mathcal{R}}(\mathbf{u})) \approx \boldsymbol{\mathcal{R}}_{ff}(\mathbf{u}_f)$
            \State Solve $\boldsymbol{\widetilde{\mathcal{R}}}_f(\mathbf{u}_f) = 0$ using few Newton iterations to get $\mathbf{u}_f^{k+1/2}$
            \State $\mathbf{u} \leftarrow \Pi_m \mathbf{u} + \Pi_f^T \mathbf{u}_f^{k+1/2}$
        \EndIf
        \State Solve $\boldsymbol{\mathcal{R}}_T(\mathbf{u}) = 0$
    \EndWhile
\EndWhile
\end{algorithmic}
\end{algorithm}

\section{Numerical Examples}
\label{sec4}
The proposed algorithm was implemented using the MATLAB Reservoir Simulation Toolbox (MRST) open-source code \citep{lie2019introduction}. Our implementation leverages the EDFM module, the automatic differentiation simulator framework, and various core utilities provided by MRST.

In this section, we present three test cases designed to evaluate the performance of the adaptive Exact Newton with Nonlinear Elimination preconditioning (EN-NE) algorithm for transport subproblems in fractured reservoirs. In each case, we compare the computational performance of EN-NE against the standard Newton method, focusing on the total number of nonlinear iterations and overall CPU time. A simple time-stepping scheme is employed, where the timestep size begins with a small initial value and progressively increases until it reaches a specified target. Relative permeabilities are modeled using a standard Corey formulation. The parameters listed in Table~\ref{tab:common-params} are common across all three numerical examples; any deviations or additional settings specific to a given case are noted accordingly within each example. For all test cases, a conservative threshold of \(\gamma = 0.25\) was adopted. In addition, as a safeguard mechanism, the preconditioner was also activated whenever the global solve goes through timestep cuts. In terms of local solves, we fix the number of iterations performed on solving the fracture subproblems to 5 for all cases. 

\renewcommand{\arraystretch}{1.2} 
\setlength{\arrayrulewidth}{0.5pt} 

\begin{table}[H]
\centering
\caption{\textit{Common parameters for all test cases.}}
\label{tab:common-params}
\begin{tabular}{@{} l c l c @{}}
\toprule
\textbf{Parameter} & \textbf{Value} & \textbf{Parameter} & \textbf{Value} \\
\midrule
Matrix porosity & 0.2 & Water compressibility &  $10^{-13}~\text{1/Pa}$ \\
Fracture porosity & 0.5 & Oil compressibility & $10^{-10}~\text{1/Pa}$ \\
Fracture aperture & 0.04 m & Fracture permeability & 1000 Darcy \\
Total simulation time, $T$ & 5 years & Target timestep, $dt$ & 30 days \\
Water viscosity & 0.001~Pa$\cdot$s & $n_w$ & 2 \\
Oil viscosity & 0.005~Pa$\cdot$s & $n_o$ & 2 \\
$\gamma$ & 0.25 & $S_{w_\text{in}}$ & 0 \\
$k^*_w$ & 1 & $k^*_o$ & 1 \\
\bottomrule
\end{tabular}
\end{table}

\subsection{Case 1}
\label{subsection4.1}
In the first test case, we consider a Cartesian grid of dimensions $120 \times 35$, featuring 48 embedded fracture planes adding 668 fracture grid cells to the computational domain. The fracture network is derived from a publicly available dataset included with the hfm module of MRST. An isotropic and homogeneous matrix permeability of 10 millidarcy is assumed. The reservoir is initially saturated with oil, and the initial pressure is set to 100 bar. To drive the flow, an injector well is placed at the bottom-right corner of the domain, injecting a total of one pore volume of water over the entire simulation time. The injector well is operated under a flow-rate control to ensure uniform injection throughout the simulation. A pressure-constrained production well is placed at the diagonally opposite corner, establishing a pressure gradient across the domain. 

\begin{figure}[H]
  \centering
  \resizebox{0.9\textwidth}{!}{%
    \begin{minipage}{\textwidth}
      \begin{minipage}[b]{0.48\textwidth}
        \centering
        \includegraphics[width=\linewidth]{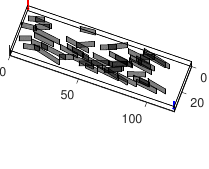}
        \subcaption{Fracture network and location of wells.}
        \label{fig:example1-setup-a}
      \end{minipage}%
      \hfill%
      \begin{minipage}[b]{0.48\textwidth}
        \centering
        \includegraphics[width=\linewidth]{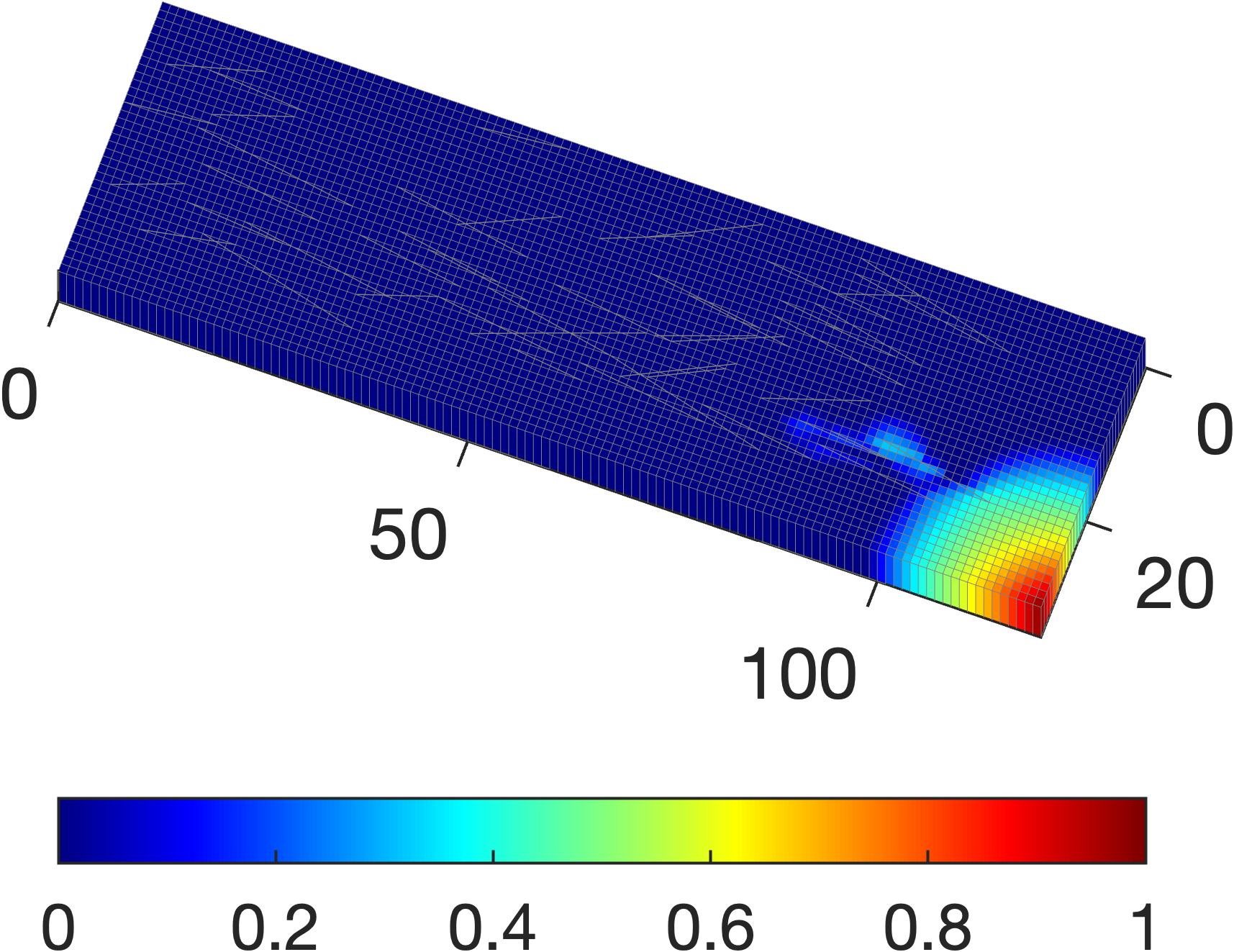}
        \subcaption{Saturation profile at \(t = 0.15T\).}
        \label{fig:example1-setup-b}
      \end{minipage}

      \vskip\baselineskip

      \begin{minipage}[b]{0.48\textwidth}
        \centering
        \includegraphics[width=\linewidth]{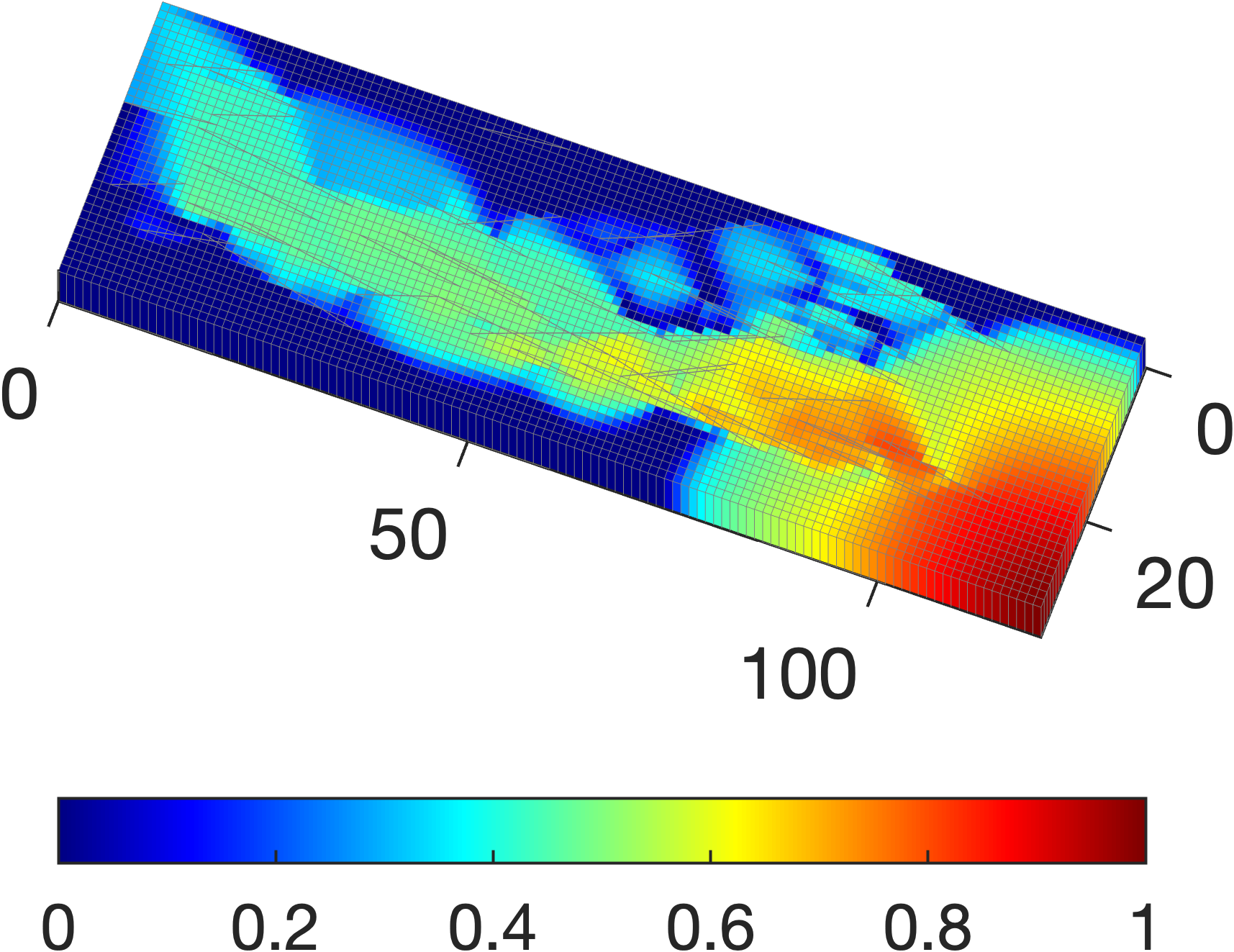}
        \subcaption{Saturation profile at \(t = 0.5T\).}
        \label{fig:example1-setup-c}
      \end{minipage}%
      \hfill%
      \begin{minipage}[b]{0.48\textwidth}
        \centering
        \includegraphics[width=\linewidth]{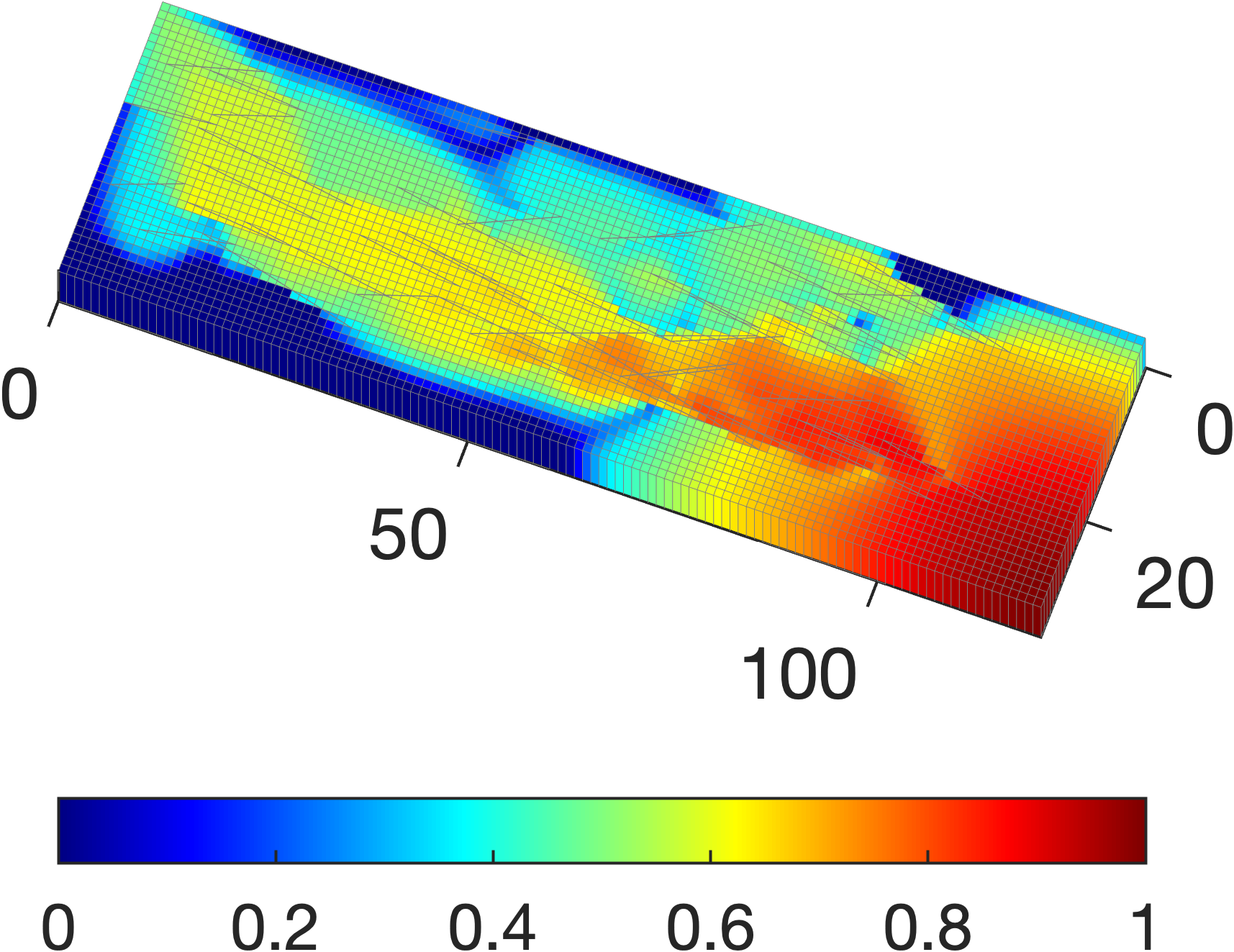}
        \subcaption{Saturation profile at \(t = T\).}
        \label{fig:example1-setup-d}
      \end{minipage}
    \end{minipage}%
  }
  \caption{Simulation setup and saturation snapshots for the first test case.}
  \label{fig:example1-setup}
\end{figure}

A schematic of the fracture network along with the injector and producer locations is provided in Figure~\ref{fig:example1-setup-a}. Additionally, Figures~\ref{fig:example1-setup}\subref{fig:example1-setup-b}, 
\ref{fig:example1-setup}\subref{fig:example1-setup-c}, and 
\ref{fig:example1-setup}\subref{fig:example1-setup-d} 
show the water saturation profiles at early, mid, and final stages of the simulation, respectively, illustrating the progression of the displacement front over time.

\begin{figure}[H]
    \centering
    \includegraphics[width=0.9\textwidth]{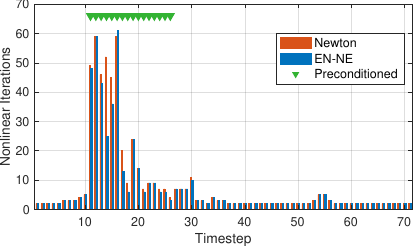}
    \caption{Nonlinear iteration count per timestep for Newton and EN-NE. Timesteps where the nonlinear preconditioner was activated are indicated by green flags.}
    \label{fig:example1-compare_nonlinear_itr}
\end{figure}

The performance of the EN-NE algorithm is evaluated for this case and benchmarked against the standard Newton method. To ensure a fair comparison, both methods utilize the default nonlinear solver settings provided in MRST. These settings include, for example, a maximum of 25 nonlinear iterations per timestep and a limit of 6 allowable timestep cuts. For further implementation details and solver configuration, the reader is referred to \citet{lie2019introduction}. Figure~\ref{fig:example1-compare_nonlinear_itr} shows the number of nonlinear iterations per timestep for both the standard Newton method and EN-NE. The green flags in the figure indicate the timesteps at which the nonlinear preconditioner was activated. Note that the preconditioner was activated only at 16 timesteps out of the total 71, demonstrating its selective and adaptive application.

\begin{figure}[H]
  \centering
  \resizebox{0.9\textwidth}{!}{%
    \begin{minipage}{\textwidth}
      \begin{minipage}[b]{0.48\textwidth}
        \centering
        \includegraphics[width=\linewidth]{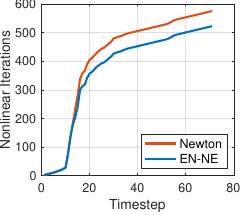}
        \subcaption{Cumulative nonlinear iterations.}
        \label{fig:example1-cum-iterations}
      \end{minipage}%
      \hfill%
      \begin{minipage}[b]{0.48\textwidth}
        \centering
        \includegraphics[width=\linewidth]{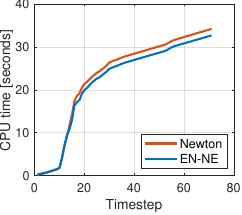}
        \subcaption{Cumulative CPU time.}
        \label{fig:example1-cum-cpu}
      \end{minipage}
    \end{minipage}%
  }
  \caption{Comparison of cumulative nonlinear iteration count and cumulative CPU time between Newton and EN-NE for the first test case.}
  \label{fig:example1-cum-performance}
\end{figure}

Despite this limited activation, EN-NE consistently required fewer nonlinear iterations compared to the standard Newton method, particularly during early timesteps where nonlinearities are typically more pronounced. Although the iteration reduction yields minute computational savings in this case, it highlights the potential of EN-NE to improve solver efficiency in more demanding scenarios. Figure~\ref{fig:example1-cum-cpu} illustrates the cumulative nonlinear iteration count and cumulative CPU time for both methods. In total, the EN-NE algorithm required 523 nonlinear iterations compared to 576 iterations for the standard Newton method. A summary of the convergence performance for all test cases is provided in Table~\ref{tab:newton-enne-results}. More pronounced computational benefits of the EN-NE method will be demonstrated in the subsequent, more complex examples.

\subsection{Case 2}
\begin{figure}[H]
  \centering
  \resizebox{0.9\textwidth}{!}{%
    \begin{minipage}{\textwidth}
      \begin{minipage}[b]{0.48\textwidth}
        \centering
        \includegraphics[width=\linewidth]{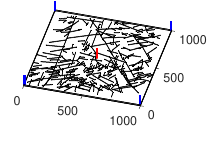}
        \subcaption{Fracture network and location of wells.}
        \label{fig:example2-setup-a}
      \end{minipage}%
      \hfill%
      \begin{minipage}[b]{0.48\textwidth}
        \centering
        \includegraphics[width=\linewidth]{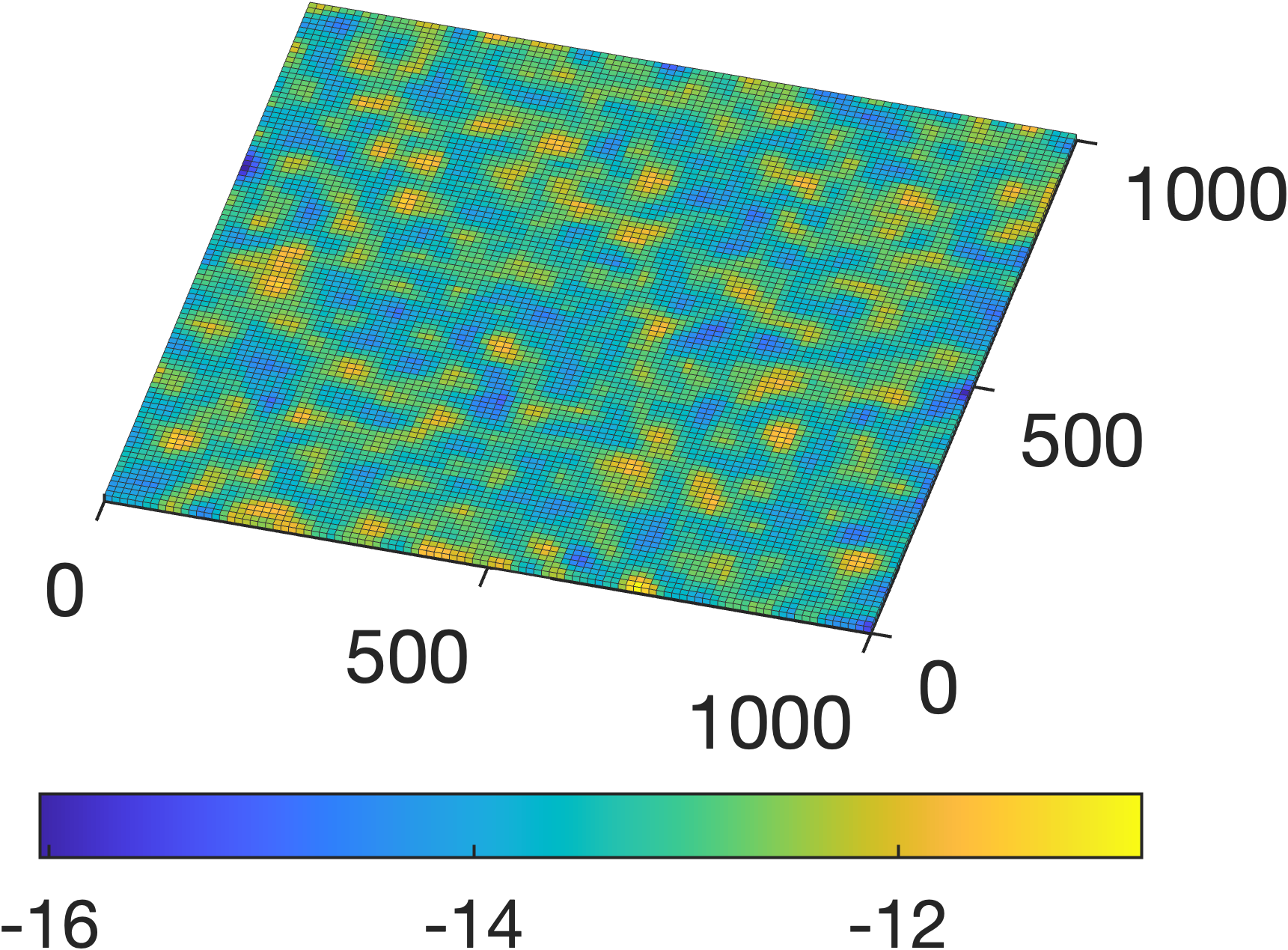}
        \subcaption{Permeability distribution, $\log(K_x)$.}
        \label{fig:example2-setup-b}
      \end{minipage}

      \vskip\baselineskip

      \begin{minipage}[b]{0.48\textwidth}
        \centering
        \includegraphics[width=\linewidth]{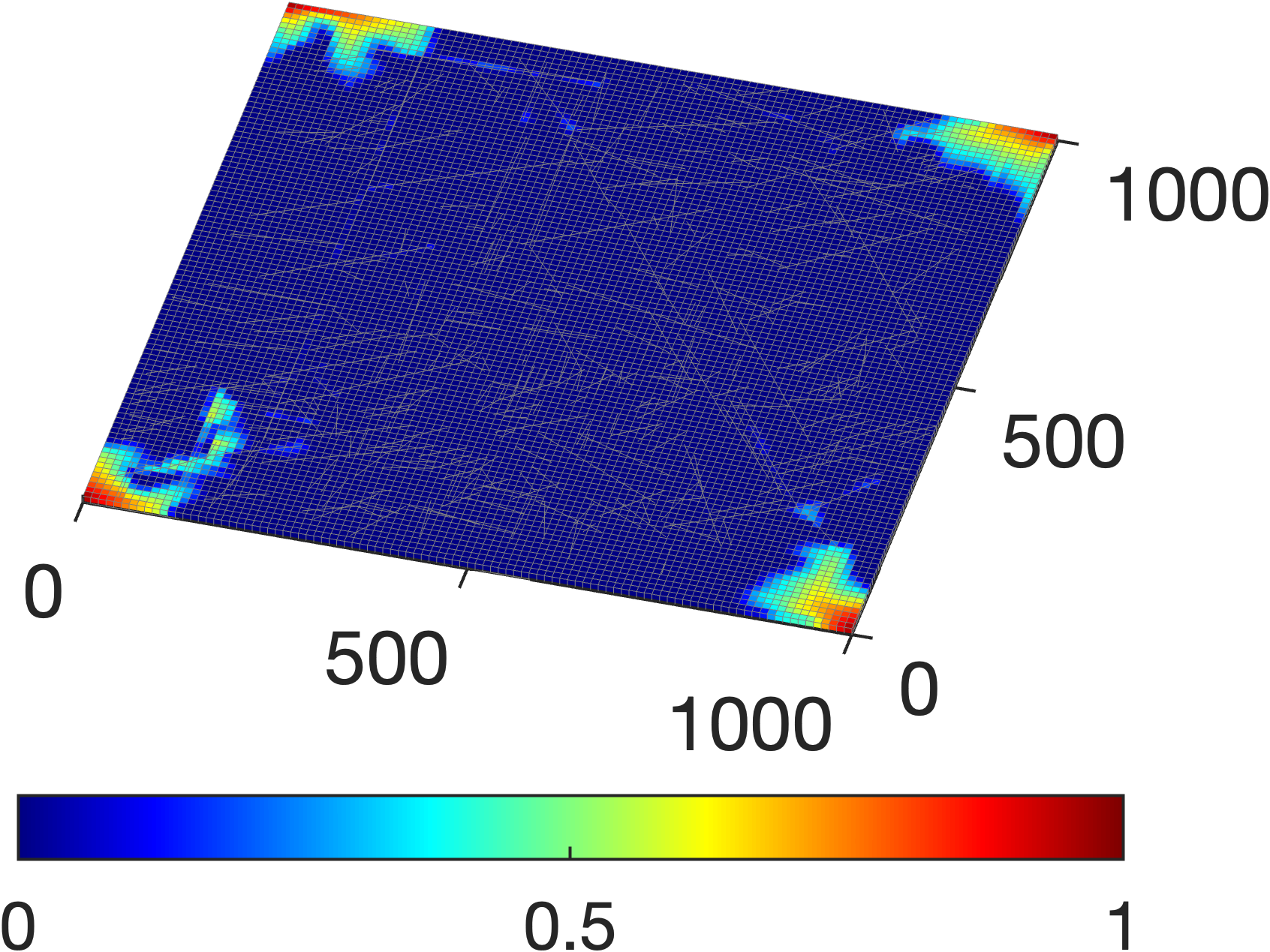}
        \subcaption{Saturation profile at \(t = 0.15T\).}
        \label{fig:example2-setup-c}
      \end{minipage}%
      \hfill%
      \begin{minipage}[b]{0.48\textwidth}
        \centering
        \includegraphics[width=\linewidth]{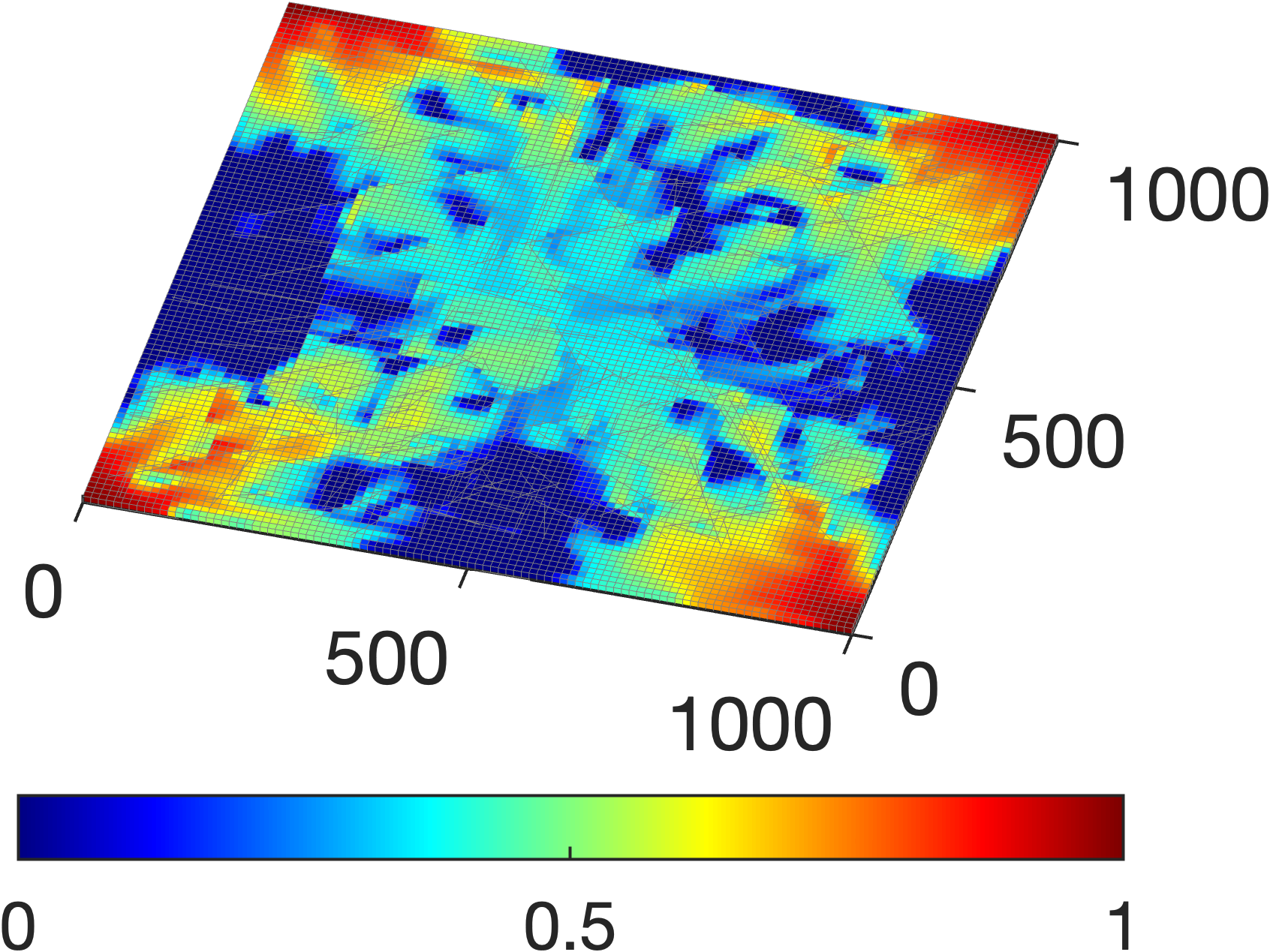}
        \subcaption{Saturation profile at \(t = 0.5T\).}
        \label{fig:example2-setup-d}
      \end{minipage}
    \end{minipage}%
  }
  \caption{Simulation setup and saturation snapshots for the second test case.}
  \label{fig:example2-setup}
\end{figure}

\begin{figure}[H]
    \centering
    \includegraphics[width=0.9\textwidth]{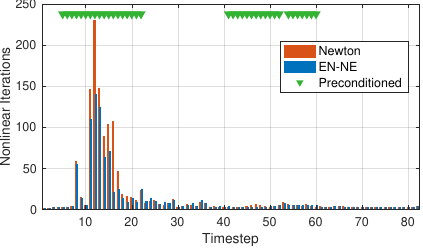}
    \caption{Nonlinear iteration count per timestep for Newton and EN-NE. Timesteps where the nonlinear preconditioner was activated are indicated by green flags.}
    \label{fig:example2-compare_nonlinear_itr}
\end{figure}

In the following test case, the convergence behavior of the EN-NE algorithm is investigated on a more complex test case. We use a 1000~m $\times$ 1000~m reservoir discretized by a $100 \times 100$ Cartesian grid. The reservoir includes a fracture network derived from outcrop data of the Jandaíra carbonate formation located in the Potiguar basin, Brazil. This dataset, originally presented in \citet{bisdom2017deterministic}, is publicly available and accessed through the MRST framework. The discretized fracture domain adds another $3529$ fracture grid cells to the computational domain. The flow is driven using a five-spot well pattern comprising four injection wells positioned at each corner of the domain and a single pressure-constrained producer well placed centrally. To deliberately induce nonlinearities through varying well controls, injection rates were set such that one pore volume of water is injected uniformly over the first 2.5 years, after which the injection rate is doubled for the remaining simulation period. Figure~\ref{fig:example2-setup-a} illustrates the fracture network along with the locations of the five wells, while Figure~\ref{fig:example2-setup}\subref{fig:example2-setup-b} provides the spatial distribution of matrix permeability. Additionally, Figures~\ref{fig:example2-setup}\subref{fig:example2-setup-c} and \ref{fig:example2-setup}\subref{fig:example2-setup-d} display the water saturation profiles at simulation times of $t=0.15T$ and $t=0.5T$, respectively.
The nonlinear iteration counts per timestep for both the standard Newton and EN-NE methods are shown in Figure~\ref{fig:example2-compare_nonlinear_itr}. Green markers identify the timesteps where the EN-NE nonlinear preconditioner was activated. Notably, in this more complex scenario, the preconditioner activation frequency increased  compared to the simpler first example. The preconditioner was activated at multiple intervals, specifically immediately following the second phase of injection around mid-simulation time.

\begin{figure}[H]
  \centering
  \resizebox{0.9\textwidth}{!}{%
    \begin{minipage}{\textwidth}
      \begin{minipage}[b]{0.48\textwidth}
        \centering
        \includegraphics[width=\linewidth]{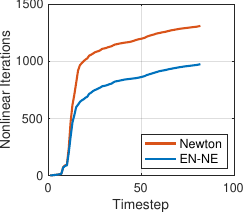}
        \subcaption{Cumulative nonlinear iterations.}
        \label{fig:example2-cum-iterations}
      \end{minipage}%
      \hfill%
      \begin{minipage}[b]{0.48\textwidth}
        \centering
        \includegraphics[width=\linewidth]{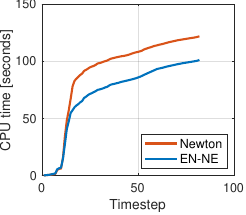}
        \subcaption{Cumulative CPU time.}
        \label{fig:example2-cum-cpu}
      \end{minipage}
    \end{minipage}%
  }
  \caption{Comparison of cumulative nonlinear iteration count and cumulative CPU time between Newton and EN-NE for the second test case.}
  \label{fig:example2-cum-performance}
\end{figure}
As can be observed, the activation of the nonlinear preconditioner substantially reduced the number of nonlinear iterations required by EN-NE relative to the standard Newton method. Figure~\ref{fig:example2-cum-performance}\subref{fig:example1-cum-iterations} depicts a notable reduction from 1310 nonlinear iterations in the standard Newton case to 975 iterations for EN-NE, representing a reduction of approximately 25\%. A significant portion of this reduction is attributed to the decrease in wasted iterations associated with timestep cuts, dropping from 450 iterations for the standard Newton method to 275 iterations for EN-NE, as detailed in Table~\ref{tab:newton-enne-results}. The corresponding CPU time plot is shown in Figure~\ref{fig:example2-cum-performance}\subref{fig:example1-cum-cpu}. For this case, EN-NE achieves a considerable CPU time improvement, reducing the total computational time from 122 seconds to 101 seconds, which corresponds to approximately 17\% reduction.
\subsection{Case 3}
\begin{figure}[H]
  \centering
  \resizebox{0.9\textwidth}{!}{%
    \begin{minipage}{\textwidth}
      \begin{minipage}[b]{0.48\textwidth}
        \centering
        \includegraphics[width=\linewidth]{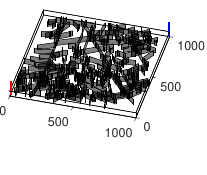}
        \subcaption{Fracture network and location of wells.}
        \label{fig:example3-setup-a}
      \end{minipage}%
      \hfill%
      \begin{minipage}[b]{0.48\textwidth}
        \centering
        \includegraphics[width=\linewidth]{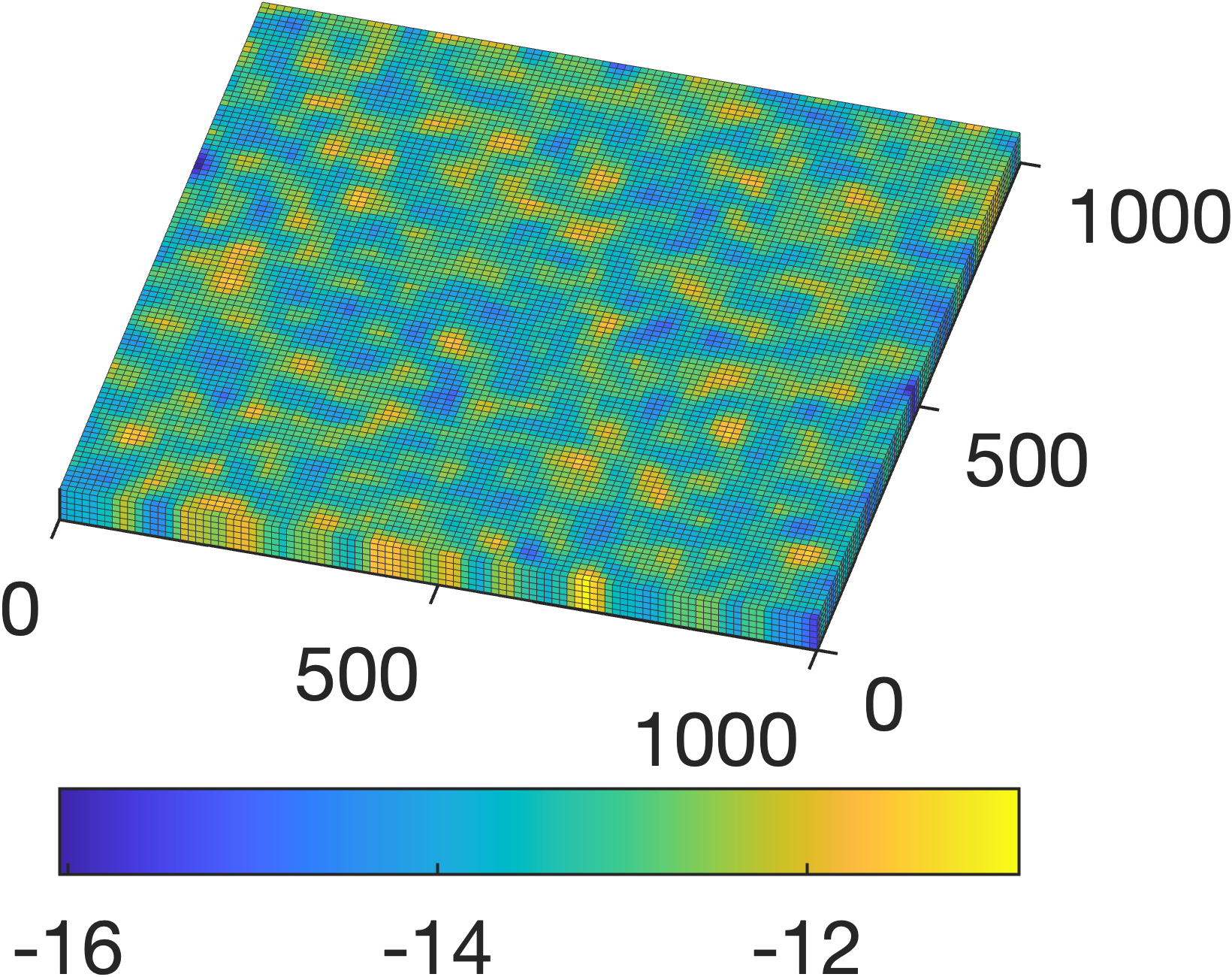}
        \subcaption{Permeability distribution, $\log(K_x)$.}
        \label{fig:example3-setup-b}
      \end{minipage}

      \vskip\baselineskip

      \begin{minipage}[b]{0.48\textwidth}
        \centering
        \includegraphics[width=\linewidth]{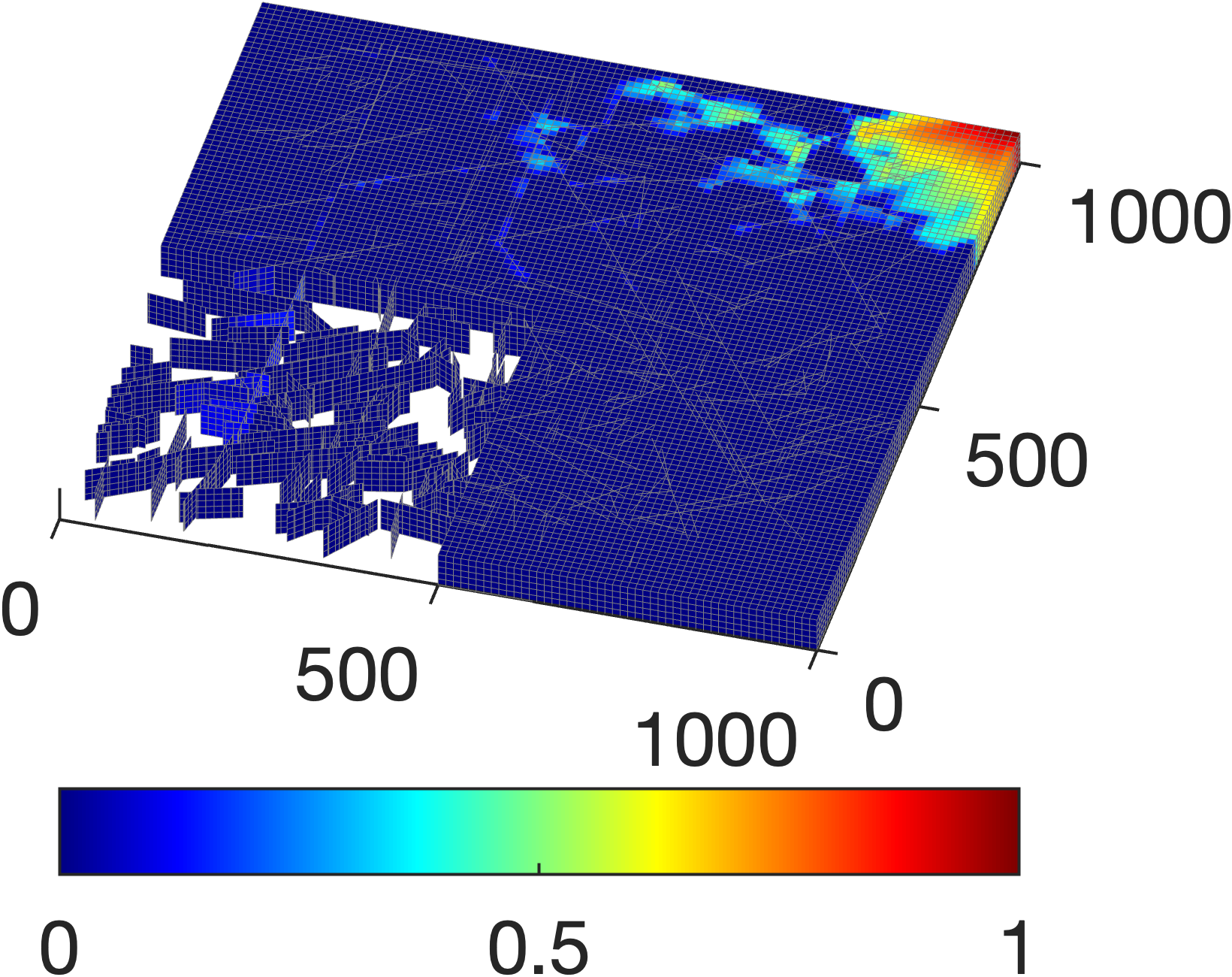}
        \subcaption{Saturation profile at \(t = 0.15T\).}
        \label{fig:example3-setup-c}
      \end{minipage}%
      \hfill%
      \begin{minipage}[b]{0.48\textwidth}
        \centering
        \includegraphics[width=\linewidth]{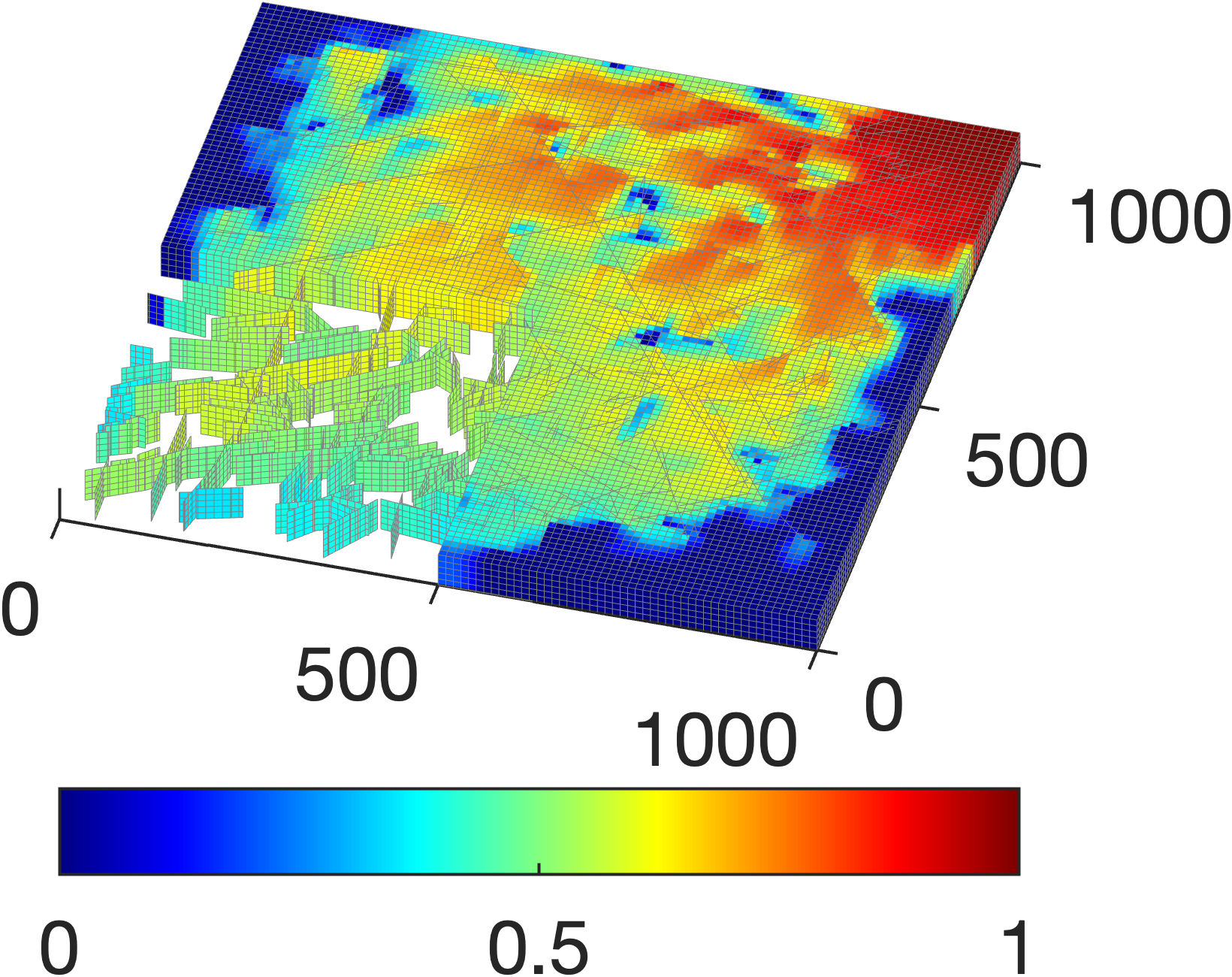}
        \subcaption{Saturation profile at \(t = T\).}
        \label{fig:example3-setup-d}
      \end{minipage}
    \end{minipage}%
  }
  \caption{Illustration of (a) the fracture network with well locations, (b) the permeability distribution, (c) the saturation profile at \(t = 0.15T\), and (d) the saturation profile at \(t = T\). A cross section of matrix cells was removed to visualize saturation profile evolving through the fracture domain.}
  \label{fig:example3-setup}
\end{figure}

To evaluate the performance of the EN-NE algorithm in a 3D fractured reservoir setup, we consider a vertically extruded version of the fracture network presented in the previous test case (Case 2). This extended fracture network is discretized along the vertical (z-axis), resulting in a computational domain composed of $100\times100\times5$ matrix grid cells and $3529\times5$ fracture grid cells. The flow scenario involves a water injection well placed at the top-right corner of the domain and a pressure-constrained producer well at the bottom-left corner. Injection occurs at a uniform and constant flow rate, delivering a total of one pore volume uniformly over a simulation period of 5 years. Figure~\ref{fig:example3-setup-a} provides a visualization of the fracture network along with the injector and producer well locations. Additionally, Figures~\ref{fig:example3-setup}\subref{fig:example3-setup-b}, \subref{fig:example3-setup-c}, and \subref{fig:example3-setup-d} illustrate, respectively, the heterogeneous distribution of matrix permeability, and the water saturation profiles at early ($t=0.15T$) and final ($t=T$) simulation run times.

\begin{figure}[H]
    \centering
    \includegraphics[width=0.9\textwidth]{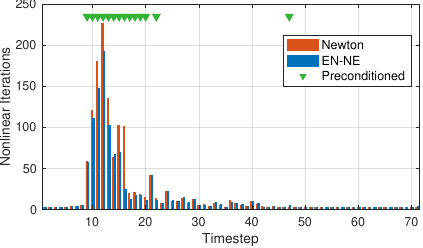}
    \caption{Nonlinear iteration count per timestep for Newton and EN-NE. Timesteps where the nonlinear preconditioner was activated are indicated by green flags.}
    \label{fig:example3-compare_nonlinear_itr}
\end{figure}

\begin{figure}[htbp]
  \centering
  \resizebox{0.9\textwidth}{!}{%
    \begin{minipage}{\textwidth}
      \begin{minipage}[b]{0.48\textwidth}
        \centering
        \includegraphics[width=\linewidth]{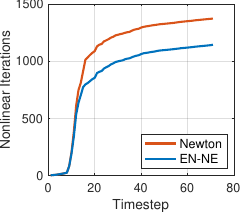}
        \subcaption{Cumulative nonlinear iterations.}
        \label{fig:example3-cum-iterations}
      \end{minipage}%
      \hfill%
      \begin{minipage}[b]{0.48\textwidth}
        \centering
        \includegraphics[width=\linewidth]{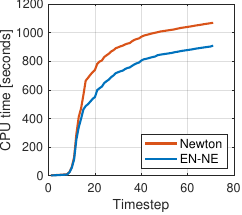}
        \subcaption{Cumulative CPU time.}
        \label{fig:example3-cum-cpu}
      \end{minipage}
    \end{minipage}%
  }
  \caption{Comparison of cumulative nonlinear iteration count and cumulative CPU time between Newton and EN-NE for the third test case.}
  \label{fig:example3-cum-performance}
\end{figure}

Figure~\ref{fig:example3-compare_nonlinear_itr} shows the nonlinear iteration counts per timestep for both the standard Newton method and EN-NE in this third example. Notably, the nonlinear preconditioner activation primarily occurred during early timesteps, a stage where localized nonlinearities were leading to bad convergence behavior for the Newton solver. The preconditioner was activated for 14 out of 71 total timesteps. Despite the limited activation of the preconditioner, EN-NE resulted in a significant reduction in the number of iterations, as illustrated in Figure~\ref{fig:example3-cum-performance}\subref{fig:example3-cum-iterations}. 

Quantitative results in Table~\ref{tab:newton-enne-results} indicate that EN-NE reduced the total nonlinear iteration count from 1347 iterations for the standard Newton solver to 1145 iterations, a 15\% reduction. A portion of this reduction is attributed to EN-NE dropping the number of wasted iterations due to timestep cuts, reducing it from 475 iterations under standard Newton to 375 iterations. Figure~\ref{fig:example3-cum-performance}\subref{fig:example3-cum-cpu} gives the cumulative CPU times—1070 seconds for Newton and 909 seconds for EN-NE, corresponding to approximately 15\% decrease in computational effort.

In all three numerical examples, our EN-NE algorithm demonstrated enhanced gains by applying the preconditioner adaptively at timesteps experiencing pronounced nonlinearities. For the physical systems considered in this work, the algorithm reduced the overall number of nonlinear iterations, including wasted iterations due to timestep cuts. The benefits of an adaptive preconditioning approach, akin to the one presented in this work, become more promising as the complexity and computational demands of the problem scale up. 

\begin{table}[H]
\centering
\caption{Comparison of total nonlinear iterations (wasted iterations on timestep cuts) and CPU time for all test cases using Newton and EN-NE.}
\label{tab:newton-enne-results}
\resizebox{\textwidth}{!}{%
\begin{tabular}{lcc@{\hskip 12pt}cc}
\toprule
& \multicolumn{2}{c}{\textbf{Newton}} & \multicolumn{2}{c}{\textbf{EN-NE}} \\
\cmidrule(lr){2-3} \cmidrule(lr){4-5}
& Iterations (wasted) & CPU time (s) & Iterations (wasted) & CPU time (s) \\
\midrule
Case 1 & 576 (150) & 34 & 523 (125) & 32 \\
Case 2 & 1310 (450) & 122 & 975 (275) & 101 \\
Case 3 & 1347 (475) & 1070 & 1145 (375) & 909 \\
\bottomrule
\end{tabular}%
}
\end{table}

\section{Concluding Remarks}

In this study, we develop an adaptive nonlinear elimination preconditioned exact Newton algorithm for the transport subproblem of two-phase flow in fractured reservoirs. Waterflooding simulations in reservoirs containing discrete fractures, modeled via EDFM, were carried out to asses the performance of the preconditioner. The adaptive NE preconditioned exact Newton algorithm consistently demonstrated favorable convergence behavior and computational efficiency compared to a standard Newton implementation. Performance gains for the studied cases include reductions in global iteration counts of up to 25\% and computational speedups of up to 17\%. 

While the current work and implementation provide a “proof-of-concept,” judicious variants of nonlinear elimination preconditioning at the algebraic and algorithmic levels are certainly possible. Potential outlooks include investigating more robust activation criteria, coupling matrix–fracture subproblem solves, employing converged (or inexact but strategically refined) local solves, and rigorously analyzing the trade-off between the additional cost of these enhanced local procedures and the attendant reduction in global iteration counts. Our ongoing research efforts include time-step optimization and investigating the use of nonlinear elimination preconditioning with inexact Newton methods for fractured reservoirs. The latter is of great interest for highly nonlinear physical problems, e.g., 3-phase and compositional modeling.

\bibliographystyle{elsarticle-harv}
\bibliography{elsarticle/bibliography}

\end{document}